\documentclass[11pt,twoside]{article}
\usepackage{a4wide}
\usepackage{amssymb}
\usepackage{amsmath}

\newtheorem{theorem}{Theorem}
\newtheorem{schmidt}{Schmidt's Theorem  \cite{WMS}}

\newtheorem{bv}{Theorem BV }

\newtheorem{mat}{Lemma \ref{slicing1}*}

\newtheorem{corollary}{Corollary}
\newtheorem{lemma}{Lemma}

\newcommand{\De}{\Delta}

\newcommand{\La}{\Lambda}

\newcommand{\nm}{{n\times m}}
\newcommand{\Rnm}{\R^{n\times m}}

\newcommand{\Inm}{\I^{n\times m}}

\renewcommand{\Bbb}[1]{\mathbb{#1}}
\newcommand{\I}{{\Bbb I}}         
\newcommand{\N}{{\Bbb N}}         
\newcommand{\R}{{\Bbb R}}         
\newcommand{\Rp}{{\Bbb R}^{+}}    
\newcommand{\Z}{{\Bbb Z}}         

\newcommand{\cA}{{\cal A}}

\newcommand{\cH}{{\cal H}}

\newcommand{\cR}{{\cal R}}

\newcommand{\cZ}{{\cal Z}}

\newcommand{\ie}{{\it i.e.}\/ }

\newcommand{\diam}{\operatorname{diam}}
\newcommand{\dist}{\operatorname{dist}}

\newcommand{\vv}[1]{{\mathbf{#1}}}

\renewcommand{\le}{\leqslant}
\renewcommand{\ge}{\geqslant}
\newcommand{\ra}{R_{\alpha}}

\begin{document}

\title{Schmidt's theorem, Hausdorff measures and Slicing}



\author{Victor Beresnevich\footnote{Research supported by EPSRC Grant R90727/01}
\\ {\small\sc York} \and  Sanju Velani\footnote{Royal Society University Research
Fellow} \\ {\small\sc York}}

\maketitle

\date{}

\centerline{{\it For Bridget Bennett on her  fortieth birthday }}

\begin{abstract} A Hausdorff measure version of
W.M. Schmidt's inhomogeneous, linear forms theorem in metric
number theory is established. The key ingredient is a `slicing'
technique motivated by  a standard result in geometric measure
theory.
In short, `slicing' together with the Mass Transference Principle
\cite{mtp}  allows us to transfer Lebesgue measure theoretic
statements for $\limsup$ sets associated with linear forms  to
Hausdorff measure theoretic statements. This extends  the approach
developed in \cite{mtp} for simultaneous approximation.
Furthermore, we establish a  new  Mass Transference Principle
which incorporates both forms of approximation. As an application
we obtain a complete metric theory for a `fully' non-linear
Diophantine problem within the linear forms setup.
\end{abstract}


\bigskip



\section{Introduction}

Fix  a vector $\vv b=(b_1,\ldots,b_m)\in\R^m$ and a non--negative,
real valued function  $\Psi:\Z^n\to\R^+ := \{ x \geq 0 : x \in \R
\} $ such that
$$\Psi(\vv a)\to0 \hspace{9mm} {\rm as \ }  \hspace{9mm} |\vv a
|:=\max(|a_1|,\ldots,|a_n|)  \, \to \, \infty \ . $$ Let
$W_{n,m}^{\vv b}(\Psi)$  be the set of  $X=(\vv x_1,\ldots,\vv
x_m)\in \Inm:=[0,1)^{\nm}$ where  $\vv x_j \in \R^n $ for  $1 \leq j
\leq m $, such that the system of inequalities
\begin{equation}\label{1}
  \|\vv a\cdot\vv x_j-b_j\| \ < \  \Psi(\vv a)  \hspace{15mm} (1
\leq j \leq m )
\end{equation}
is satisfied for infinitely many $\vv a\in\Z^n $. Here and
throughout $\vv x\cdot\vv y=x_1y_1+\dots+x_ny_n$ is the standard
inner product of two vectors $\vv x$,$\vv y \in \R^n$ and $\|x\|$
is the distance from $x \in \R$ to the nearest integer. The
following is a geometric interpretation of the set $W_{n,m}^{\vv
b}(\Psi)$  and brings to the forefront its $\limsup$ nature. For
vectors $\vv a\in\Z^n $, $\vv b\in\R^m$ and $\vv p   \in \Z^m$,
consider the $(n-1)m$--dimensional plane $ R_{{\vv a},{\vv
p}}^{\vv b}$ given by
\begin{equation}
\label{resset}
 R_{{\vv a},{\vv p}}^{\vv b} := \{ X \in \Rnm: \
\vv a\cdot\vv x_j-b_j = p_j \ \  (1 \leq j \leq m ) \, \} \ \ .
\end{equation}
Thus, $R_{{\vv a},{\vv p}}^{\vv b}$ is the product of the
$n$--dimensional hyperplanes given by $$ R_{{\vv a},{p_j}}^{b_j}
:= \{ \vv x \in \R^n: \ \vv a\cdot\vv x - b_j = p_j  \} \qquad 1
\leq j \leq m  \ . $$ For $\delta \geq 0 $, let $\Delta( R_{{\vv
a},{\vv p}}^{\vv b}, \delta ) $ denote the $\delta$--neighborhood
of $R_{{\vv a},{\vv p}}^{\vv b}$; i.e.
the product of the
$\delta$-neighborhoods (with respect to the Euclidean norm) of the
 hyperplanes $ R_{{\vv a},{p_j}}^{b_j} $. Note that
when $n=1$, the $\delta$--neighborhood $\Delta( R_{{\vv a},{\vv
p}}^{\vv b}, \delta ) $ is simply a ball of radius $\delta$ in the
supremum norm.   
It is easily verified, that
$$ X \in W_{n,m}^{\vv b}(\Psi) \hspace{10mm}  {\rm if \ and \ only \ if  }  \hspace{10mm}
X \in \Delta\Big( R_{{\vv a},{\vv p}}^{\vv b},
\textstyle{\frac{\Psi(\vv a)}{\sqrt{{\vv a}.{\vv a}}} }  \Big)
\cap \Inm
$$ for infinitely many vectors $\vv a\in\Z^n$ and $ \vv p \in
\Z^m$.

%
%


\subsection{The Lebesgue measure theory: Schmidt's Theorem}
The following key result   provides a beautiful and simple
criteria for the  `size' of the set $W_{n,m}^{\vv b}(\Psi)$
expressed in terms of $\nm$--dimensional Lebesgue measure $|\  \
|_{\nm}$.  The theorem is due to W.M. Schmidt  and shows that
$|W_{n,m}^{\vv b}(\Psi)|_{\nm}$  satisfies an elegant `zero--one'
law.

\begin{schmidt}
Let $\Psi$ be as above and $n+m>2$. Then
$$
|W_{n,m}^{\vv b}(\Psi)|_{\nm}=\left\{
\begin{array}{rl}
0 & {\rm if} \;\;\;
\displaystyle \sum_{\vv a\in\Z^n}\Psi(\vv a)^m<\infty\,\\[4ex]
1 & {\rm if} \;\;\; \displaystyle \sum_{\vv a\in\Z^n}\Psi(\vv
a)^m=\infty\,
\end{array}
\right.. $$
\end{schmidt}

In fact, Schmidt considers a more general setup for which he obtains
a quantitative result. The case that $m+n=2$, corresponding to
$n=m=1$, is  naturally excluded since  the statement is known to be
false -- the Duffin--Schaeffer conjecture \cite{mtp,Spr79} provides
the appropriate `expected' statement.

In order to appreciate  the true  significance  of Schmidt's
Theorem it is well worth mentioning  a few special cases which in
their own right represent landmarks  within the classical theory
of metric Diophantine approximation.
\begin{description}
\item[~\hspace*{3ex} Khintchine's Theorem  \cite{Kh} \ : ] $n=1$ and ${\vv b}= 0 $ with $\Psi$
monotonic. ~\vspace*{-1ex}

~\hspace*{32ex} (simultaneous, homogeneous approximation)
\item[~\hspace*{3ex} Groshev's Theorem  \cite{groshev} \ \hspace*{3ex}: ] $n > 1$ and ${\vv b}= 0
$ with $\Psi$ monotonic\footnote{When $n>1$, the function  $\Psi$
is in general multi-variable.  To say that it is monotonic simply
means that $\Psi(\vv a):= \psi(|\vv a|)$ for some  monotonic
$\psi$. }.~\vspace*{-1ex}


~\hspace*{32ex} (linear forms, homogeneous approximation)
\item[~\hspace*{3ex} Gallagher's Theorem \cite{gal} \ \hspace*{1.2ex}: ] $n=1$ and $m \geq
2$. ~\vspace*{-1ex}

~\hspace*{32ex} (simultaneous, inhomogeneous approximation)
\end{description}

\noindent Note that under the condition that $\Psi$ is monotonic,
the results of Khintchine and Groshev already give the homogeneous
version (${\vv b}= 0 $)  of Schmidt's Theorem without the condition
that $m+n > 2 $  -- Khintchine's Theorem covers the case $n=m=1$.
Generalizing the Khintchine--Groshev statement to the inhomogeneous
case  and entirely removing the monotonicity condition when $m +n
> 2 $ is by no means a trivial feat. As with Schmidt, Gallagher considers
a more general setup for which he obtains a quantitative result.

\subsection{The general metric theory}
Let $f$ be a dimension function and let $\cH^f$ denote the Hausdorff
$f$--measure -- see \S\ref{HM}.  In short, our aim is to provide a
complete metric theory for the set $W_{n,m}^{\vv b}(\Psi)$. The
following result achieves this goal in that it provides a  simple
criteria for the `size' of the set $W_{n,m}^{\vv b}(\Psi)$ expressed
in terms of the general measure $\cH^f$.

\begin{theorem}\label{t1}
Let $\Psi$  be as above and $n+m>2$. Let $f$ be a dimension
function such that $r^{-nm}f(r)$ is monotonic. Furthermore, assume
that $g : r \to r^{-(n-1)m} f(r)$ is a dimension function. Then
$$
\cH^f(W_{n,m}^{\vv b}(\Psi))=\left\{
\begin{array}{cl}
0& {\rm if}  \qquad\displaystyle \sum_{\vv
a\in\Z^n\smallsetminus\{\vv 0\}} \
g\!\left(\dfrac{\Psi(\vv a)}{|\vv a|}\right)  \times \ |\vv a|^m \ < \ \infty\, \\[4ex]
\cH^f(\Inm)  & {\rm if}  \qquad\displaystyle\sum_{\vv
a\in\Z^n\smallsetminus\{\vv 0\}} \ g\!\left(\dfrac{\Psi(\vv
a)}{|\vv a|}\right)  \times \ |\vv a|^m \ = \ \infty\,
\end{array}
\right..
$$
\end{theorem}

Notice that in the case $\cH^f$ is $\nm$--dimensional Lebesgue
measure $|\  \ |_{\nm}$, the theorem reduces to Schmidt's Theorem.
As with the  Lebesgue theory, the convergence part of the above
theorem is relatively straightforward if not trivial -- see
\S\ref{t1conv}. The main substance is the divergent part.  For
this, our particular strategy is straightforward enough. We
establish the following:
\begin{theorem}\label{t2}
$
 ~ \ \  {\it Schmidt's \ Theorem      \  (divergent \ part) }  \ \ \
\Longrightarrow \ \ \  {\it Theorem \  \ref{t1}  \ (divergent \
part) } \ . $
\end{theorem}


In \cite{mtp}, this strategy has recently  been successfully
implemented to establish the simultaneous version of Theorem
\ref{t2}; namely

\begin{bv}\label{bv}
$ {\it Gallagher's \ Theorem      \  (divergent \ part) } \
\Longrightarrow  \    {\it Theorem \  \ref{t1}  \ (divergent \
part) }
\\ ~\hspace*{66ex} with \   n=1 \  and  \  m \geq 2  \ .
$
\end{bv}

Recall, that Schmidt's Theorem reduces to Gallagher's Theorem in
the case of simultaneous approximation ($n=1$). To be absolutely
precise, in \cite[\S6.2]{mtp}  we  only  consider the homogeneous
case 
of Theorem BV. However, given the method of proof adopted in
\cite{mtp}  no extra obstacles appear in establishing the
inhomogeneous version above. Indeed, the proof is essentially a
simple application of the Mass Transference Principle (see
\S\ref{secmtp}) and for this it is irrelevant whether we start
with a homogeneous or inhomogeneous divergent statement of
Gallagher's Theorem. The only relevant aspect is that when $n=1$,
the set $W_{n,m}^{\vv b}(\Psi)$ is a limsup set naturally defined
in terms  of a sequence of balls -- the Mass Transference
Principle then does the rest!    This is no longer the case when
$n>1$ and so Theorem \ref{t2}  is not simply a consequence of the
approach developed in \cite{mtp};  namely the Mass Transference
Principle.

The key aspect of this paper is the introduction of a `slicing'
technique to the theory of metric Diophantine approximation; in
particular to the linear forms aspect of the theory.  The
technique is motivated by a relatively standard result in
geometric measure theory -- see \S\ref{secslice}. The upshot is
that `slicing' together with the Mass Transference Principle
yields Theorem \ref{t2} -- the `hard' part of Theorem \ref{t1}.

\vspace{1ex}

{\em Remark. \ }  In all previous contributions towards the general
metric theory, such as the pioneering work of Jarn\'{\i}k \cite{Ja}
the function $\Psi$ is assumed to be monotonic. For further details
and references the reader is refereed to \cite[Sections 1.1 \&
12.1]{BDV03}.

\vspace{1ex}

Before moving on, it is useful to say a little
concerning the condition imposed on $g$  in Theorem \ref{t1};
namely that since $g$ is assumed to be a  dimension function we
have that $g(r) \to 0$ as $r \to 0 $. For the sake of clarity and
ease of discussion, put $f: r \to r^s$ ($s>0$) in Theorem
~\ref{t1}. Then, Theorem~\ref{t1} reduces to the following
$s$--dimensional Hausdorff measure 
statement which in its own right is of significant importance
since it characterizes the Hausdorff dimension of the set
$W_{n,m}^{\vv b}(\Psi)$ as the exponent of convergence of a
certain `$s$--volume' sum.

\begin{corollary}\label{cor1}
Let $\Psi$  be as above and $n+m>2$. Let $\delta>0$ and
$s:=(n-1)m+\delta$. Then
$$
\cH^s(W_{n,m}^{\vv b}(\Psi))=\left\{
\begin{array}{cl}
0& {\rm if}  \qquad\displaystyle \sum_{\vv
a\in\Z^n\smallsetminus\{\vv 0\}} \
\Psi(\vv a)^\delta  \ |\vv a|^{m-\delta} \ < \ \infty\, \\[4ex]
\cH^s(\Inm)  & {\rm if}  \qquad\displaystyle\sum_{\vv
a\in\Z^n\smallsetminus\{\vv 0\}}\ \Psi(\vv a)^\delta  \ |\vv
a|^{m-\delta} \ = \ \infty\,
\end{array}
\right..
$$
\end{corollary}

\noindent It follows from the definition of Hausdorff dimension
(see \S\ref{HM})  that if  for some  $\delta>0$ 
the sum in the corollary diverges, then
$$
\dim W_{n,m}^{\vv b}(\Psi) \ = \ \inf \left\{ s:
\textstyle{\sum_{\vv a\in\Z^n\smallsetminus\{\vv 0\}} } \ \Psi(\vv
a)^{s-(n-1)m} \ |\vv a|^{nm- s} \ < \ \infty \right\}  \ .
$$

We suspect that the condition on $s$ imposed in  Corollary
\ref{cor1}, namely that $s$ is strictly greater than  $(n-1)m$,
cannot be relaxed. Briefly, if $\delta = 0$ so that $s=(n-1)m$, the
sum in Corollary~\ref{cor1} diverges irrespective of  $\Psi$. Now
the `approximating'
hyperplanes 
as defined by (\ref{resset}) are themselves of dimension $s$ and
indeed  of positive $\cH^s$ measure. Thus, it is highly likely
that for rapidly decreasing functions $\Psi$ the $\cH^s$--measure
theoretic structure of $W_{n,m}^{\vv b}(\Psi)$ is purely dependent
on  the arithmetic properties of the approximating hyperplanes. In
view of this, for appropriate $\Psi$ and  $\vv b\in\R^m$ one might
expect that $\cH^s (W_{n,m}^{\vv b}(\Psi)) $ is finite and
possibly even zero rather than $\cH^s(\Inm)$ which is infinite.


%
%
%
%
%
%
%

\section{Preliminaries}

\subsection{Hausdorff measures \label{HM}}

In this section we give a brief account of Hausdorff measures. For
further details see \cite{MAT}. A {\em dimension function} $f \, :
\, \R^+ \to \R^+ $ is a continuous, non-decreasing function such
that $f(r)\to 0$ as $r\to 0 \, $.
The Hausdorff
$f$--measure with respect to the dimension function $f$ will be
denoted throughout by ${\cal H}^{f}$ and is defined as follows.
Suppose $F$ is  a  subset of $\R^k$. For $\rho
> 0$, a countable collection $ \left\{B_{i} \right\} $ of balls in
$\R^k$ with radius  $r(B_i) \leq \rho $ for each $i$ such that $F
\subset \bigcup_{i} B_{i} $ is called a {\em $ \rho $-cover for
$F$}.
 For a dimension
function $f$ define $$
 {\cal H}^{f}_{\rho} (F) \, = \, \inf \ \sum_{i} f(r(B_i)),
$$
where the infimum is taken over all $\rho$-covers of $F$. The {\it
Hausdorff $f$--measure} $ {\cal H}^{f} (F)$ of $F$ with respect to
the dimension function $f$ is defined by   $$ {\cal H}^{f} (F) :=
\lim_{ \rho \rightarrow 0} {\cal H}^{f}_{\rho} (F) \; = \;
\sup_{\rho > 0 } {\cal H}^{f}_{\rho} (F) \; . $$

A simple consequence of the definition of $ {\cal H}^f $ is the
following useful fact.

\begin{lemma}
If $ \, f$ and $g$ are two dimension functions such that the ratio
$f(r)/g(r) \to 0 $ as $ r \to 0 $, then ${\cal H}^{f} (F) =0 $
whenever ${\cal H}^{g} (F) < \infty $. \label{dimfunlemma}
\end{lemma}

In the case that  $f(r) = r^s$ ($s > 0$), the measure $ \cH^f $ is
the usual $s$--dimensional Hausdorff measure $\cH^s $ and the
Hausdorff dimension $\dim F$ of a set $F$ is defined by $$ \dim \, F
\, := \, \inf \left\{ s : \mathcal{ H}^{s} (F) =0 \right\} = \sup
\left\{ s : \mathcal{ H}^{s} (F) = \infty \right\} . $$ In
particular when $s$ is an integer, $\cH^s$ is
comparable\footnote{The symbols $\ll$ and $\gg$ will be used to
indicate an inequality with an unspecified positive constant. If $a
\ll b $ and $a \gg b $ we write $a \asymp b $, and say that the
quantities $a$ and $b$ are comparable.  } to the $s$--dimensional
Lebesgue measure. Actually, $\cH^s$ is a constant multiple of the
$s$--dimensional Lebesgue measure.

\subsection{The  Mass Transference Principle \label{secmtp}}

Given a  dimension function $f$ and a ball $B=B(x,r)$ in $\R^m$,
we define another ball
\begin{equation}\label{e:001}
\textstyle B^f:=B(x,f(r)^{1/m}) \ .
\end{equation}
When $f(x)=x^s$ for some $s>0$ we also adopt the notation $B^s$, \ie
$ B^s:=B^{(x\mapsto x^s)}. $ It is readily verified that
\begin{equation}\label{e:001a}
B^m=B.
\end{equation}

Given a sequence of balls $B_i$, $i=1,2,3,\ldots$, as usual its
limsup set is
$$
\limsup_{i\to\infty}B_i:=\bigcap_{j=1}^\infty\ \bigcup_{i\ge j}B_i \
.
$$
By definition,  $\limsup_{i\to\infty}B_i$ is precisely the set of
points in $\R^m$ which   lie in infinitely many balls $B_i$.

\vspace*{1ex}

The following  Mass Transference Principle allows us to  transfer
Lebesgue measure theoretic statements for $\limsup$ subsets of
$\R^m$ to Hausdorff measure theoretic statements.

\begin{lemma}[Mass Transference Principle]\label{thm3}
Let $\{B_i\}_{i\in\N}$ be a sequence of balls in $\R^m$ with
$r(B_i)\to 0$ as $i\to\infty$. Let $f$ be a dimension function such
that $r^{-m}f(r)$ is monotonic and let $\Omega$ be a ball in $\R^m
$. Suppose for any ball $B$ in $\Omega$
$$
\cH^m\big(B \, \cap \, \limsup_{i\to\infty}B^f_i{}\,\big) \ = \
\cH^m(B) \ .
$$
Then, for any ball $B$ in $\Omega$
$$
\cH^f\big(B \, \cap \, \limsup_{i\to\infty}B^m_i\,\big)\ = \
\cH^f(B) \ .
$$
\end{lemma}

%
%
%

With $\Omega = \R^m$, the lemma is precisely   Theorem~2 in
\cite{mtp}. It is easily seen that this implies the above modified
statement which is better suited for the particular applications we
have in mind.


\subsection{The `slicing' lemma \label{secslice} }

The `slicing' lemma below  is the crucial new ingredient and
together with the `slicing' technique (\S\ref{secst})  makes it
possible to reduce Theorem~\ref{t2} to an $m$-dimensional problem by
slicing the original set $W_{n,m}^{\vv b}(\Psi)$ into a family of
subsets lying on parallel $m$-dimensional planes. The `slicing'
technique is motivated by the  `slicing' lemma.

 In the following
$V$ will be a linear subspace of $\R^k$ and $V^\perp$ will denote
the linear subspace of $\R^k$ orthogonal to $V$. Also,
$V+a:=\{v+a:v\in V\}$ for $ a \in V^\perp$.


\begin{lemma} \label{slicing1}
Let $l,k\in\N$ such that $l \leq k $ and $f$ and $g:r\mapsto
r^{-l}f(r)$ be dimension functions. Let $A\subset\R^k$ be a Borel
set with $\cH^f(A)<\infty$. Then for any $(k-l)$-dimensional
linear subspace $V$ of\/ $\R^k$,
$$
\cH^{g}(A\cap(V+a))<\infty  \ \  \  \text{for $\cH^{l}$-almost all
$a\in V^\perp $}.
$$
\end{lemma}

\vspace{2ex}

When $f: r \to r^s$, the  lemma constitutes the first part of
Theorem~10.10 in \cite{MAT}. The proof given there can be easily
modified to yield the more general statement above. Nevertheless,
given the importance of the lemma  and for the sake of
completeness we have included the proof of Lemma \ref{slicing1} as
an appendix.

\vspace{2ex}



Trivially, Lemma~\ref{slicing1} implies the following:

\begin{lemma}[Slicing lemma]\label{slicing2}
Let $\, l,k\in\N$ such that $l \leq k $ and $f$ and $g:r\mapsto
r^{-l} f(r)$ be dimension functions. Let $A\subset\R^k$ be a Borel
set and $V$ be an $(k-l)$-dimensional linear subspace of $\R^k$.
If for a subset $S$ of $V^\perp$ of positive $\cH^{l}$-measure
$$
\cH^{g}(A\cap(V+b))=\infty\text{ \ \ \ \  for all \ \  }b\in S \,
,
$$
then $\cH^f(A)=\infty$.
\end{lemma}

\subsection{Additional assumption in Schmidt's theorem (divergent part) \label{add}}

 Let $n \geq 2 $ as otherwise the substance of this
section becomes trivial. For each $i\in\{1,\dots,n\}$ define the
subset $\cZ_i$ of $\Z^n\smallsetminus\{\vv0\}$ to consist of vectors
$\vv a\in \Z^n$ such that $|\vv a|=|a_i|$. Assume that
$$
\sum_{\vv a\in\Z^n}\Psi(\vv a)^m=\infty.
$$
Now
$$
\infty=\sum_{\vv a\in\Z^n\smallsetminus\{\vv0\}}\Psi(\vv a)^m\le
\sum_{i=1}^n\sum_{\vv a\in\cZ_i}\Psi(\vv a)^m.
$$
Therefore there is an index $i\in\{1,\dots,n\}$ such that
\begin{equation}\label{v}
\sum_{\vv a\in\cZ_i}\Psi(\vv a)^m=\infty.
\end{equation}
Define
$$
\Psi_i(\vv a)= \left\{
\begin{array}{cl}
 \Psi(\vv a), & \text{if }\vv a\in\cZ_i,\\[2ex]
 0 , & \text{if }\vv a\not\in\cZ_i.
\end{array}\right.
$$

Hence, $|W_{n,m}^{\vv b}(\Psi_i)|_{\nm}=1$ by Schmidt's Theorem.
Trivially, this implies that for almost all $X\in\Inm$
\begin{equation}\label{e:004}
\max_{1\le j\le m}\|\vv a\cdot\vv x_j-b_j\|\ < \ \Psi(\vv a)
\end{equation}
for infinitely many $\vv a\in\cZ_i$. There is no loss of
generality in assuming that (\ref{v}) is satisfied with $i=1$ as
otherwise we can apply a permutation of variables (columns in $X$)
under which Schmidt's theorem is clearly invariant. Thus, when
considering the divergent part of Schmidt's theorem we can assume
that
\begin{equation*} \label{e:A083}
 \Psi (\vv a )  \, = \, 0  \  \ \ \ \ \forall  \ \ \vv  a \in \Z^n
 \hspace{5ex}  {\rm with } \hspace{5ex}  |\vv a | \neq |a_1|  \ .
\end{equation*}

\section{Proof of Theorem~\ref{t1}}

\subsection{The case of convergence \label{t1conv}}

 We are given that
$$
\sum_{\vv a\in\Z^n\smallsetminus\{\vv 0\}} \  g\!
\left(\dfrac{\Psi(\vv a)}{|\vv a|}\right) \times |\vv a|^m  \ <  \
\infty \ .
$$

The convergent part of Theorem \ref{t1} follows on  using standard
covering arguments. For each $N \in \N $, it is easily verified that
$$
W_{n,m}^{\vv b}(\Psi) \ \subset \ \bigcup_{\substack{ {\vv a}  \in
\Z^n: |{\vv a}| \geq N \\  \Psi(\vv a) > 0 }} \ \  \bigcup_{{\vv p}
\in \Z^m } \ \Delta\Big( R_{{\vv a},{\vv p}}^{\vv b},
\textstyle{\frac{\Psi(\vv a)}{|{\vv a }|} } \Big) \cap \Inm \ . $$
 Note that there is no loss of generality in assuming that
 $\Psi(\vv a)>0$ in the above union, for otherwise (\ref{1}) has no solutions
$X$  and the integer vector $\vv a$ makes no contribution to
$W_{n,m}^{\vv b}(\Psi)$.

 Next notice that for any fixed ${\vv a}  \in
\Z^n\smallsetminus\{\vv0\}$ with $\Psi(\vv a)
> 0 $ and ${\vv p} \in \Z^m $, it is possible to cover
$$\Delta\Big( R_{{\vv a},{\vv p}}^{\vv b}, \textstyle{\frac{\Psi(\vv a)}{|{\vv a
}|} } \Big) \cap \Inm  $$ by a collection ${\cal C }_{{\vv a},{\vv
p}}^{\vv b} $ of balls of common radius $ \frac{\Psi(\vv a)}{|{\vv
a }|} $ such that $$ \# \, {\cal C }_{{\vv a},{\vv p}}^{\vv b} \ll
\Big(\textstyle{\frac{|{\vv a }|}{\Psi(\vv a)} }\Big)^{(n-1)m} \ \
.
$$ Also, for a fixed ${\vv a} \in \Z^n\smallsetminus\{\vv0\}$
$$ \# \, \left\{ {\vv p} \in \Z^m  : \, \Delta\Big( R_{{\vv a},{\vv p}}^{\vv b},
\textstyle{\frac{\Psi(\vv a)}{|{\vv a }|} } \Big) \cap \Inm \ \neq
\emptyset  \ \right\} \ \ll \ |{\vv a}|^m  \ . $$ Finally, note
that since $\Psi({\vv a})  \to 0 $ as $ |{ \vv a } | \to \infty $,
we have that for all $N$ sufficiently large
$$
\frac{\Psi(\vv a)}{|{\vv a }|}  \ \leq \frac{1}{N} \  \ . $$

It now follows from the definition of ${\cal H}^f $ that for $N$
sufficiently large,
\begin{eqnarray*}
{\cal H}^f_{\rho := \frac{1}{N} } \Big(W_{n,m}^{\vv b}(\Psi)\Big)
& \ll & \sum_{\substack{ {\vv a}  \in \Z^n: |{\vv a}| \geq N }}
f\!\left(\dfrac{\Psi(\vv a)}{|\vv a|}\right) \times
\left(\dfrac{\Psi(\vv a)}{|\vv a|}\right)^{-(n-1)m}  \times \ |\vv
a|^m \\ & & \\ & := & \sum_{\substack{ {\vv a}  \in \Z^n: |{\vv
a}| \geq N }} g\!\left(\dfrac{\Psi(\vv a)}{|\vv a|}\right) \times
\ |\vv a|^m \  \ \ \to \  \ 0 \hspace{4mm} {\rm as \ }
\hspace{3mm} N \to \infty \ .
\end{eqnarray*}
Thus, ${\cal H}^f (W_{n,m}^{\vv b}(\Psi)) = 0 $ as required.
\hfill $\Box$

\subsection{The case of divergence (Theorem \ref{t2}): the `slicing' technique \label{secst} }
 \emph{Throughout we assume that $n\geq 2$.
Theorem BV covers the  $n=1$ case.} We are given that
\begin{equation} \label{divv}
\sum_{\vv a\in\Z^n\smallsetminus\{\vv 0\}} \  g\!
\left(\dfrac{\Psi(\vv a)}{|\vv a|}\right) \times |\vv a|^m \ = \
\infty \ .
\end{equation}

We start by  considering the case that $r^{-mn}f(r) \to L $ as $r
\to 0 $ and $L$ is finite. If $L=0$, then Lemma \ref{dimfunlemma}
implies that $\cH^f(\Inm)=0$ and since $ W_{n,m}^{\vv b}(\Psi)
\subset \Inm $ the result follows. If $L \neq 0 $, then $\cH^f$ is
comparable to $\cH^{mn}$ (in fact, $\cH^f = L \, \cH^{mn}$). In
turn, $\cH^{mn}$ is comparable to $\nm$--dimensional Lebesgue
measure and so the required statement follows on showing that
$|W_{n,m}^{\vv b}(\Psi)|_{\nm}  =   |\Inm|_{\nm}  $. Well, this is
a simple consequence of Schmidt's theorem since the sum appearing
in (\ref{divv}) is comparable to $\sum_{\vv a\in\Z^n} \Psi(\vv
a)^m $.

In view of the above discussion, we can assume without loss of
generality that
\begin{equation} \label{maininf}
r^{-mn}f(r) \ \to \ \infty \hspace{6mm} {\rm as } \hspace{6mm}
r\to0 \ \ . \end{equation}
 Indeed, it is this situation that
constitutes the main substance of Theorems \ref{t1} and \ref{t2}.
Trivially, (\ref{maininf}) together with Lemma \ref{dimfunlemma}
implies that
$$ \cH^f(\Inm) \ = \ \infty \  .  $$  Similarly,  since 
(\ref{maininf}) is equivalent to the statement that $r^{-m}g(r)
\to \infty  $ as $ r \to 0$, we have that $\cH^g(B)=\infty$ for
any $m$--dimensional ball $B$.

To proceed, we set
$$
\widetilde\Psi(\vv a)^m  \ := \  g \! \left(\dfrac{\Psi(\vv
a)}{|\vv a|}\right)\times |\vv a|^m \ .
$$
In view of (\ref{divv}), it follows  that
$$
\sum_{\vv a\in\Z^n}\widetilde\Psi(\vv a)^m=\infty \ ,
$$
 and  Schmidt's Theorem (divergent part)   implies  that
$$
|W_{n,m}^{\vv b}(\widetilde\Psi)|_{\nm}=1 \ .
$$
The goal is to show that this implies that
$$\cH^f(W_{n,m}^{\vv b}(\Psi)) \ =  \ \infty  \  ; $$
i.e. to establish Theorem \ref{t2} under the condition imposed by
(\ref{maininf}). Recall, that the conclusion of  Theorem \ref{t2}
is  precisely the divergent part of Theorem \ref{t1}.


In view of the discussion in  \S\ref{add}, we can assume without
loss of generality that
\begin{equation} \label{e:083}
 \widetilde\Psi (\vv a )  \, = \, 0  \  \ \ \ \ \forall  \ \ \vv  a \in \Z^n
 \hspace{5ex}  {\rm with } \hspace{5ex}  |\vv a | \neq |a_1|  \ .
\end{equation}
Now, let
$$
V=\{(\vv x_1,\dots,\vv x_m) \  :\   x_{j,i}=0 \  \ \forall\
j={1,\dots,m}\; ; i={2,\dots,n},\ \},
$$
where $\vv x_j=(x_{j,1},\dots,x_{j,n})$. Thus, $V$ is an
$m$--dimensional subspace of $\Rnm$. By Fubini's theorem, there is
a subset $S\subset \I^{m(n-1)}\subset V^\perp$  with $|S|_{m(n-1)}
= 1 $ such that for every $X_0\in S$ the set $W_{n,m}^{\vv
b}(\widetilde\Psi)$ has full $m$-dimensional Lebesgue measure in $
(V+{X_0})\cap\Inm $; i.e.
\begin{equation}\label{e:080}
|(V+X_0)\cap W_{n,m}^{\vv b }(\widetilde\Psi) |_m \ = \ 1  \ .
\end{equation}
For every $(\vv p,\vv a)=(p_{1},\dots p_{m}; a_1,\dots,a_n)\in
\Z^{m} \times \Z^{n}$ define the set $\sigma(\vv p,\vv a)$ to
consist of $X\in\Inm$ such that
\begin{equation*}
\max_{1\le j\le m}|\vv a\cdot\vv x_j-b_j+p_{j}|<\widetilde\Psi(\vv
a)  \ .
\end{equation*}
In view of condition (\ref{e:083}) imposed on $\widetilde\Psi$,
the set $\sigma(\vv p,\vv a)$ is empty whenever $|\vv a | \neq
|a_1|$. We therefore assume that $|\vv a | =  |a_1|$ throughout
the rest of the proof.  Then
\begin{equation}\label{e:006}
\sigma(\vv p,\vv a)\cap (V+{X_0})
\end{equation}
is the product of $m$ intervals of length $2\widetilde\Psi(\vv
a)/|a_1|$. Indeed, for all $j=1,\dots,m$ and $i=2,\dots,n$ we have
that  $x_{j,i}$ are fixed and defined by $X_0$ for all  points in
this set. That is to say that the only coordinates that may vary
are $x_{j,1}$. Therefore, the set (\ref{e:006}) is defined by the
system
\begin{equation*}
\max_{1\le j\le
m}|a_1x_{j,1}+(a_2x_{j,2}+\dots+a_nx_{j,n}+p_{j}-b_j)|<\widetilde\Psi(\vv
a)
\end{equation*}
or equivalently
\begin{equation}\label{e:008}
\max_{1\le j\le
m}\left|x_{j,1}-\frac{b_j-(a_2x_{j,2}+\dots+a_nx_{j,n}+p_{j})}{a_1}\right|
<\frac{\widetilde\Psi(\vv a)}{|a_1|} \ .
\end{equation}
The pathological situation of $a_1=0$ is excluded by the
conditions $|\vv a|=|a_1|$ and $\vv a\not=\vv0$.  On identifying
$\Inm$ with the $\nm$-dimensional torus it is easily seen that
every inequality in (\ref{e:008}) defines an interval of length
$2\widetilde\Psi(\vv a)/|a_1|$. Thus (\ref{e:006}) defines  a ball
of radius $\widetilde\Psi(\vv a)/|\vv a|$. Such balls form a
sequence $(A_i)_{i\in\N}$. Therefore
$$
\limsup_{i\to\infty}A_i=(V+X_0)\cap W_{n,m}^{\vv b
}(\widetilde\Psi) \ .
$$
Hence, in view of (\ref{e:080}) we have that
$|\limsup_{i\to\infty}A_i|_m=1$. This implies that   for any ball
$B\subset (V+X_0)\cap \Inm$
$$
\cH^m(\limsup_{i\to\infty}A_i\cap B)=\cH^m(B)\,.
$$

\noindent For each ball $A_i$ define  the ball $B_i$ with the same
centre and radius $\Psi(\vv a)/|\vv a|$. Then, by definition
$B_i^{g}=A_i$ -- see (\ref{e:001}). It follows that
\begin{equation}
\label{e:081} \limsup_{i\to\infty} B_i \ \subset \ (V+X_0)\cap
W_{n,m}^{\vv b }(\Psi)  \ .
\end{equation}
Also,   $r^{-m} g(r) = r^{-mn} f(r) $   is monotonic by assumption.
Thus, on applying  the Mass Transference Principle with $ \Omega =
(V+X_0)\cap \Inm$, we obtain that for any ball $B$ in $\Omega$
$$
\cH^{g}(\limsup_{i\to\infty}B_i^m\cap B)=\cH^{g}(B)=\infty \ .
$$
Recall, that $\cH^{g}(B)=\infty$ is a consequence of (\ref{maininf})
and Lemma \ref{dimfunlemma}.  Hence, in view of (\ref{e:081}) and
the fact that $B_i^m := B_i$ we have that for every  $X_0\in S$
$$
\cH^{g}((V+X_0)\cap W_{n,m}^{\vv b}(\Psi))=\infty \ .
$$
Recall that $\cH^{m(n-1)}(S) > 0$ and so by the Slicing lemma,
$$\cH^f(W_{n,m}^{\vv b}(\Psi))=\infty  \ .  $$
This completes the proof of  Theorem \ref{t2} and therefore  the
divergent part of Theorem \ref{t1}.

\hfill $\Box$

\section{A Mass Transference Principle for linear forms}

The Mass Transference Principle 
deals with $\limsup$ sets which are defined as a sequence of
balls. However, we have seen that together with the `slicing'
technique introduced in this paper  we are able to deal with
$\limsup$ sets defined as a sequence of neighborhoods of
`approximating' planes -- at least within the context of Schmidt's
Theorem. In short, the aim of this section is to develop a single
framework which enables us to combine the Mass Transference
Principle and `slicing' into a single statement. The main result
(Theorem \ref{t3} below) should be viewed as a generalization to
the linear forms setup of the Mass Transference Principle
developed in \cite{mtp} for simultaneous approximation. As
applications, we deduce Theorem \ref{t2} (which constitutes the
main substance of Theorem \ref{t1})  as a simple corollary and
more strikingly, we obtain a complete metric theory for a `fully'
non-linear Diophantine problem -- see Theorem \ref{t4} of
\S\ref{fnl}. 


\subsection{A general framework for approximating  by planes  \label{gf}}

Throughout  $k, m \geq 1 $ and $ l\ge0$ are integers such that
$k=m+l$.
%
%
%
%
Let $\cR=(\ra )_{\alpha \in J}$ be a family of planes in $\R^k$ of
common dimension $l$ indexed by an infinite countable set $J$. For
every $\alpha\in J$ and $\delta\geq 0$ define
$$ \Delta(R_\alpha,\delta) := \{\vv x \in \R^k: \dist(\vv
x,R_\alpha) < \delta\} \ . $$ Thus $\Delta(R_\alpha,\delta)$ is
simply the $\delta$--neighborhood of the $l$--dimensional plane
$R_\alpha$. Note that by definition, $ \Delta(R_\alpha,\delta) =
\emptyset $ if $\delta =0$. Next, let
$$\Upsilon : J \to \R^+ : \alpha\mapsto
\Upsilon(\alpha):=\Upsilon_\alpha$$ be a non-negative, real valued
function on $J$. In order to avoid pathological situations within
our framework, we assume that for every $\epsilon
>0$ the set $\{\alpha\in J:\Upsilon_\alpha>\epsilon \}$ is finite.
This condition implies that $\Upsilon_\alpha  \to 0 $ as $\alpha$
runs through $J$. We now consider the following `$\limsup$' set,
$$ \La(\Upsilon)=\{\vv x\in\R^k:\vv
x\in\Delta(R_\alpha,\Upsilon_\alpha)\ \mbox{for\ infinitely\ many\
}\alpha\in J\} \ . $$ Note that in view of the conditions imposed
on $k,l$ and $m$ we have that $l < k$. Thus the dimension of the
`approximating' planes $ R_\alpha$ is strictly less than that of
the ambient space $\R^k$. The situation when $l =k $ is of little
interest and has therefore been naturally omitted.

\subsection{The main result}

\begin{theorem}\label{t3}
Let $\cR$ and $\Upsilon$ as above be given. Let $V$ be a linear
subspace of\/ $\R^k$ such that $\dim V=m=\mathrm{codim}\,\cR$ and

~ \qquad  \qquad $(i)$\quad \ $V \ \cap \  R_\alpha \ \neq \
\emptyset $ \quad for all $ \ \alpha\in J \ $,

~ \qquad  \qquad $(ii)$\quad $\sup_{\alpha\in J}\diam( \,
V\cap\De(R_\alpha,1) \, ) \ < \ \infty \ $ .

\noindent Let $f$ and $g : r \to g(r):= r^{-l} \, f(r)$ be
dimension functions such that $r^{-k}f(r)$ is monotonic and let
$\Omega $ be a ball in $\R^k$. Suppose for any ball $B$ in
$\Omega$ $$ \cH^k \big( \, B  \cap \La \big(g(\Upsilon)^{\frac1m}
\big) \,  \big) \, = \, \cH^k(B)
$$ Then
$$ \cH^f \big( \, B \cap \Lambda(\Upsilon) \, \big) \, = \, \cH^f(B) \ .  $$
\end{theorem}

\vspace{2ex}

{\it Remark} : Conditions  (i) and (ii) are not particularly
restrictive. When $l=0$,  so that  $\cR$ is a collection of points
in $\R^k $, conditions (i) and (ii) are trivially satisfied and
Theorem \ref{t3} simply reduces to the Mass Transference Principle
of \S\ref{secmtp}. When $l \geq 1 $, so that $\cR$ is a collection
of $l$--dimensional planes  in $\R^k $, condition (i) excludes
planes $\ra$  parallel to $V$ and condition (ii) simply means that
the angle at which  $R_\alpha$ `hits'  $V$ is bounded away from zero
by a fixed constant independent of $\alpha \in J$. This in turn
implies that each plane in $\cR$ intersects $V$ at exactly one
point.


\vspace{2ex}


\subsection{Theorem \ref{t3} $\Longrightarrow $ Theorem \ref{t2} }
 With reference to our general framework, let $k = m
\times n$. Hence, $l= m(n-1)$. Furthermore, let $ J := \{ (\vv a,
\vv p, \vv b) \in   \Z^n\setminus\{\vv 0\}  \times \Z^m \times
\{\vv b \} : |\vv a | = |a_1| \} $ where $\vv b $ is a fixed
vector in $\R^m$, $\alpha := (\vv a, \vv p, \vv b) \in J $,
$R_\alpha := R_{{\vv a},{\vv p}}^{\vv b}$ where the latter is
given by (\ref{resset}) and $ \Upsilon_\alpha :=
\textstyle{\frac{\Psi(\vv a)}{\sqrt{{\vv a}.{\vv a}}} } $. Then,
$$ W_{n,m}^{\vv b}(\Psi)  \  \supset \
\widetilde{W}_{n,m}^{\vv b}(\Psi)  \ :=  \ \La(\Upsilon)  \, \cap \, \Inm \ \
. $$ In view of \S\ref{add}, it suffices to establish Theorem
\ref{t2} for the set $\widetilde{W}_{n,m}^{\vv b}(\Psi)$.  As in
\S\ref{secst},  let
$$
V : =\{(\vv x_1,\dots,\vv x_m) \  :\   x_{j,i}=0 \  \ \forall\
j={1,\dots,m}\; ; i={2,\dots,n},\ \}  \ ,
$$
where $\vv x_j=(x_{j,1},\dots,x_{j,n})$. Thus, $V$ is an
$m$--dimensional subspace of $\Rnm$ and  it is easily verified that
conditions (i) and (ii) of Theorem \ref{t3} are satisfied. Theorem
\ref{t2} now follows on applying Theorem \ref{t3} with $\Omega =
\Inm $. Note that we can actually deduce the following stronger
`local' statement.  For any ball $B$ in $ \Inm $,
$$
\cH^f(B \cap W_{n,m}^{\vv b}(\Psi)) \, = \, \cH^f(B)   \hspace{6ex}
{\rm if}  \hspace{3ex} \qquad\displaystyle\sum_{\vv
a\in\Z^n\smallsetminus\{\vv 0\}} \ g\!\left(\dfrac{\Psi(\vv a)}{|\vv
a|}\right)  \times \ |\vv a|^m \ = \ \infty  \ \ .
$$

\subsection{Preliminaries}

Before embarking on the proof of Theorem \ref{t3}, we derive some
crucial facts from conditions (i) and (ii) imposed in the
statement of the theorem. We also state a  `shrinking' lemma which
will be required in the proof of Theorem \ref{t3}.

\subsubsection{Crucial consequences of conditions (i) and (ii)}

Let $l \geq 1 $ as otherwise the substance of this section is
irrelevant. Thus, $\cR$ is a family of `genuine' planes and not
points. In view of the remark immediately after the statement of
Theorem \ref{t3}, we have that for every $\alpha\in J$ there is
the unique point $c_\alpha$ given by  $ V\cap R_\alpha$.  Clearly,
the ball $B(c_\alpha,r)$ in $\R^k$ is contained in the
$r$--neighborhood of $R_\alpha$; i.e. $B(c_\alpha,r)  \subset
\De(R_\alpha,r)$. Hence,
$$B(c_\alpha,r) \, \cap \,  V \ \subset \  \De(R_\alpha,r) \, \cap \,  V   \ . $$
It follows that the diameter of $\De(R_\alpha,r)\cap V$ is at
least $2r$ -- the diameter of the ball $B(c_\alpha,r)\cap V$.
On the other hand, condition (ii) implies that the diameter of
$\De(R_\alpha,r)\cap V$ is bounded above by a constant $C>0$ times
$r$ (uniformly in $\alpha$). Indeed, with $C$  equal to the
supremum in the left hand side of condition (ii) we have that
$$
\De(R_\alpha,1) \, \cap \,  V \ \subset \  B(c_\alpha,C) \, \cap
\,  V \ .
$$
Since $R_\alpha$ and $V$ are planes, the set $\De(R_\alpha,r)\cap
V$ is simply the set $\De(R_\alpha,1)\cap V$  scaled by the factor
$r$ -- shrunk or expanded depending on whether $r$ is less than or
greater  than  one. Similarly,  $B(c_\alpha,Cr)\cap V$  is
$B(c_\alpha,C)\cap V$ scaled by the factor $r$.   The upshot of
this, is that

\begin{equation}\label{vb1}
B(c_\alpha,r)\cap V\ \subset\ \De(R_\alpha,r)\cap V\ \subset\
B(c_\alpha,Cr)\cap V\qquad \text{for any } \ r>0\ .
\end{equation}

Finally, we observe that since $R_\alpha$ and $V$ are planes, the
inclusions given by  (\ref{vb1}) remains valid  if $V$ is replaced
by any parallel hyperplane. Formally, for any $r>0$ and any $\vv
x_0\in\R^k$

\begin{equation}\label{vb2}
B(c_{\alpha,\vv x_0},r)\cap (V+\vv x_0)\ \subset\
\De(R_\alpha,r)\cap (V+\vv x_0)\ \subset\ B(c_{\alpha,\vv
x_0},Cr)\cap (V+\vv x_0)\ ,
\end{equation}
where $c_{\alpha,\vv x_0}$ is the  unique point given by
$R_\alpha\cap (V+\vv x_0)$.

\subsubsection{The shrinking lemma}

Given a ball $B$ and a positive constant $\delta<1$, let $\delta
B$  denote the ball $B$ shrunk by the factor $\delta$. The
following result formally states that the measure of $\limsup$
sets arising from a sequence of balls in $\R^k$ is not effected by
shrinking the balls  by a constant factor.

\begin{lemma} \label{shrlem}
Let $B_i$ be a sequence of balls in $\R^k$ such that
$\limsup_{i\to\infty} B_i$ has full measure in an open subset $U$
of $\R^k$. Let  $\delta<1$ be a positive constant.  Then, the set
$\limsup_{i\to\infty} \delta B_i$ has full measure in $U$.
\end{lemma}

The lemma is  a simple consequence of Lemma~6 in \cite{mtp}.


\subsection{Proof of Theorem~\ref{t3}}

\emph{Without loss of generality, we assume that $l \geq 1 $. The
case when $l=0$ corresponds  to the Mass Transference Principle of
\S\ref{secmtp}.  }

\vspace{1ex}

The proof of Theorem \ref{t3}
 follows the basic strategy as  the proof of Theorem \ref{t2}
(see \S\ref{secst}); i.e. that of combining the  Mass Transference
Principle (Lemma \ref{thm3}) and the Slicing lemma (Lemma
\ref{slicing2}) in an appropriate manner.  In view of this we
shall give a sketch proof and leave the details to the reader.

As in the proof of Theorem \ref{t2}, we can assume without loss of
generality that
\begin{equation} \label{maininf1}
r^{-k}f(r) \ \to \ \infty \hspace{6mm} {\rm as } \hspace{6mm}
r\to0 \ \ . \end{equation} Indeed, it is this situation that
constitutes the main substance of Theorem \ref{t3}. Recall, that
(\ref{maininf1}) together with Lemma \ref{dimfunlemma} implies
that $\cH^f(B)=\infty$ for any $k$--dimensional ball $B$ and that
$\cH^g(B)=\infty$ for any $m$--dimensional ball $B$. For the sake
of clarity we introduce the following notation. Let $V$ be as in
the statement of the theorem. For a subset $A$ of $\R^k$ and $\vv
x_0$ in $ V^\perp$ let
$$
A'_{\vv x_0}  \; := \; A \, \cap \, (V+\vv x_0)   \ \ . $$ By
definition,
$$
\La'_{\vv x_0}(\Upsilon)=\{\vv x\in\R^k:\vv x\in\Delta'_{\vv x_0}
(R_\alpha,\Upsilon_\alpha)\ \mbox{for\ infinitely\ many\
}\alpha\in J\} \ .
$$

Fix a ball $D$ in $\Omega$. The aim is to show that
$$ \cH^f(D \cap \Lambda(\Upsilon)) \, = \, \infty  \ .  $$
We are given that \begin{equation} \label{sv1}
 \cH^k\big(D \cap
\La\big(g(\Upsilon)^{\frac1m}\big)\big) \, = \, \cH^k(D) \ .
\end{equation} Now let $ D^* := \{ \vv x_0 \in V^\perp : D'_{\vv x_0} \neq \emptyset
\} $.  Then, (\ref{sv1}) together with Fubini's theorem implies
the existence of  a  set $S\subset D ^* \subset V^\perp$ with
$|S|_{l} =  |D^*|_{l} $ such that for every $ {\vv x_0 } \in S$
\begin{equation}\label{vb8}
\cH^m\big(D'_{\vv x_0 } \; \cap \; \La'_{\vv x_0
}\big(g(\Upsilon)^{\frac1m}\big) \big) \ = \  \cH^m(D'_{\vv x_0 })
\ .
\end{equation}
In view of (\ref{vb2}),  we have that
\begin{equation}\label{vb4}
\limsup_{\alpha \in J} \  B'_{\vv x_0 } \big(c_{\alpha}^*,
g(\Upsilon_\alpha)^{\frac1m}\big) \ \subset \  \La'_{\vv x_0
}\big(g(\Upsilon)^{\frac1m}\big) \ \subset  \ \limsup_{\alpha \in J}
\ B'_{\vv x_0 }\big(c_{\alpha}^*,Cg(\Upsilon_\alpha)^{\frac1m}\big)
\ \ .
\end{equation}
For each $\alpha \in J$,  the ball $ B'_{\vv x_0 } (c_{\alpha}^*,r)
$ is by definition a subset of   $V +{\vv x_0 }$ with centre $
c_{\alpha}^* := R_\alpha\cap (V+\vv x_0)$.  It follows via
(\ref{vb8}) and (\ref{vb4}),  that
\begin{equation}\label{vb7}
\cH^m\big(\, D'_{\vv x_0 }  \cap \; \textstyle{\limsup_{\alpha \in
J}} \ B'_{\vv x_0
}\big(c_{\alpha}^*,Cg(\Upsilon_\alpha)^{\frac1m}\big) \  \big) \ = \
\cH^m(D'_{\vv x_0 }) \ .
\end{equation}
As a consequence of the shrinking lemma (Lemma \ref{shrlem}), we
can put $C=1$ in (\ref{vb7}); i.e.
\begin{equation}\label{vb6}
\cH^m\big( \, D'_{\vv x_0 }  \cap \; \textstyle{\limsup_{\alpha \in
J}}  \ B'_{\vv x_0 }\big(c_{\alpha}^*,
g(\Upsilon_\alpha)^{\frac1m}\big) \ \big) \ = \ \cH^m(D'_{\vv x_0 })
\ .
\end{equation}
Now  for any ball $B$ in $D'_{\vv x_0 }$, (\ref{vb6}) implies that
$$
\cH^m\big( \, B \cap \; \textstyle{\limsup_{\alpha \in J}}  \
B'_{\vv x_0 }\big(c_{\alpha}^*, g(\Upsilon_\alpha)^{\frac1m}\big) \
\big) \ = \ \cH^m(B ) \ .
$$
On applying the Mass Transference Principle  with $\Omega =
D'_{\vv x_0 }$, we obtain that
\begin{equation}\label{vb5}
\cH^g\big( \, D'_{\vv x_0 } \,  \cap \, \textstyle{\limsup_{\alpha
\in J}} \ B'_{\vv x_0 }(c_{\alpha}^*, \Upsilon_\alpha) \ \big) \; =
\; \cH^g (D'_{\vv x_0 })  \; = \;   \infty \ .
\end{equation}
In view of (\ref{vb2}),  we have that
\begin{equation*}\label{vb3}
\limsup_{\alpha \in J} \; B'_{\vv x_0 }(c_{\alpha}^*,
\Upsilon_\alpha) \ \subset \ \La'_{\vv x_0 }(\Upsilon) \ \subset \
\limsup_{\alpha \in J} \; B'_{\vv x_0 }(c_{\alpha}^*,
C\Upsilon_\alpha)   \ .
\end{equation*}
This together with (\ref{vb5}), implies that for every ${\vv x_0 }
\in S $
$$
\cH^g\big( \, D'_{\vv x_0 } \,  \cap \, \La'_{\vv x_0 }(\Upsilon)
\ \big) \; = \;   \infty \ .
$$
On applying the Slicing lemma, we obtain that $ \cH^f(D \cap
\Lambda(\Upsilon))  =  \infty $ as desired.

\hfill $\Box$

\subsection{`Fully' non-linear Diophantine problems \label{fnl}}

Schmidt's theorem underpins the metric theory of non-linear
Diophantine approximation --  the  integer points $\vv a$
associated with the definition of $W_{n,m}^{\vv b}(\Psi)$ can be
restricted to lie in a subset $\cA$ of $\Z^n$ which is completely
free of any linear structure. Indeed, one  simply sets $\Psi$ to
be zero for points $\vv a$ outside of $\cA$ so that the points
$\vv a$ that make any contribution to $W_{n,m}^{\vv b}(\Psi)$ lie
 only in $\cA$. However, Schmidt's theorem is not non-linear in the
full sense, since the integer variable $\vv p$, implicit in the
symbol $\|\cdot\|$ is a linear term. Theorems~\ref{t1} is
therefore of the same nature; i.e. it provides a complete metric
theory of non-linear Diophantine approximation but fails to be
fully non-linear. Sprind\v{z}uk, in his 1979 monograph
\cite{Spr79} writes: `As of now, no metric theory of (fully)
non-linear Diophantine approximation has been constructed. The
working out of such a theory is a very topical problem.' Since
then, substantial progress has been made within the one
dimensional setting  -- the numerator and denominator of the
rational approximates $a/p$ are restricted to sets of number
theoretic interest such as primes (see, for example \cite[Chapter
6]{har}  for the Lebesgue measure theory and \cite[\S12.5]{BDV03}
for the complete metric theory). There has also been some progress
within the simultaneous setting \cite{jones}. To our knowledge,
there has been no progress what so ever within the linear forms
setting.  We now demonstrate the power of Theorem~\ref{t3} -- it
naturally allows us to consider fully non-linear problems; in
particular within the linear forms setting.

%
%
%

A  natural source of fully non-linear Diophantine problems is the
theory of partial differential equations (PDE's). 
The following is a concrete  example of a fully non-linear problem
arising in such a manner -- it is related to the solubility of the
two-dimensional inhomogeneous wave equation (see \cite{BDKL} for
details). Given a vector $\vv a=(a_1,a_2)\in\Z^2$,  let $\vv
a^2:=(a_1^2,a_2^2)$. Let $\psi:\Rp\to\Rp$ be a non-negative,
 real valued function and consider the set
$$
S_2(\psi):=\{\vv x\in\I^2:|\,\vv a^2\cdot\vv x-p^2|<\psi(|\vv a|)
\text{ for infinitely many }(\vv a,p)\in\Z^2\times\Z\ \}\ .
$$
Naturally, the problem is to determine a complete metric theory for
$S_2(\psi)$. Clearly, this is a fully non-linear problem since the
coefficients of the `approximating planes'  are restricted to
perfect squares.


\medskip

In \cite{BDKL}, the following criteria for the `size' of the set
$S_2(\psi)$ expressed in terms of $2$--dimensional Lebesgue
measure  $|\  \ |_{2}$ is established.

\noindent{\bf Theorem BDKL } {\it Let $\psi:\Rp\to\Rp$  be a
monotonic function such that $\lim_{h\to\infty}\psi(h)=0$. Then $$
|S_2(\psi)|_2  =\left\{
\begin{array}{rl}
0&  {\rm if}  \qquad 
\sum_{h=1}^\infty  \ \psi(h) \, < \, \infty\,\\[2ex]
1& {\rm if}  \qquad 
\sum_{h=1}^\infty \ \psi(h) \, = \, \infty\,
\end{array}
\right.. $$ }

\medskip

In view of our general framework and  Theorem~\ref{t3}, we are
able to give a complete measure theoretic description of the set
$S_2(\psi)$.


\medskip

\begin{theorem} \label{t4} Let $\psi:\Rp\to\Rp$  be a monotonic function such that
$\lim_{h\to\infty}\psi(h)=0$. Let $f$ be a dimension function such
that $r^{-2}f(r)$ is monotonic. Furthermore, assume that $g : r
\to r^{-1} f(r)$ is a dimension function. Then
$$
\cH^f(S_2(\psi))=\left\{
\begin{array}{cl}
0& {\rm if}  \qquad \displaystyle \sum_{h=1}^\infty  \
g\!\left(\dfrac{\psi(h)}{h^2}\right)  \times \ h^2 \ < \ \infty\, \\[4ex]
\cH^f(\I^2)  & {\rm if}  \qquad \displaystyle \sum_{h=1}^\infty \
g\!\left(\dfrac{\psi(h)}{h^2}\right) \times \ h^2 \ = \ \infty\,
\end{array}
\right..
$$
\end{theorem}

%

With $f:r \to r^s  \;  (s > 0) $, the theorem reduces the the
following $s$--dimensional Hausdorff measure statement. Naturally,
 it coincides with  Theorem BKDL when $s=2$.

\begin{corollary} \label{c4} Let $\psi:\Rp\to\Rp$  be a monotonic function such that
$\lim_{h\to\infty}\psi(h)=0$. For   $ 1 < s \leq 2$, we have that
$$ \cH^s(S_2(\psi))=\left\{
\begin{array}{cl}
0& {\rm if}  \qquad 
\sum_{h=1}^\infty \  \psi(h)^{s-1}h^{4-2s} \, < \, \infty\,\\[2ex]
\cH^s(\I^2) &  {\rm if}  \qquad
\sum_{h=1}^\infty \ \psi(h)^{s-1}h^{4-2s} \, = \, \infty\,
\end{array}
\right..
$$
\end{corollary}

Consider the case  $\psi : r \to r^{-\tau}  $  ($\tau > 0$) and
write $ S_2(\tau)$  for $S_2(\psi) $. For $\tau > 1$,  the above
corollary not only implies that
$$
\dim S_2(\tau) \ = \ \textstyle{ \frac{5+ \tau}{2 + \tau } }  \ ,
$$
but that $\cH^s(S_2(\tau))$ is infinite at the critical exponent
$s= \dim S_2(\tau)$.

\subsubsection{Proof of Theorem \ref{t4} } We start be rewriting the set $S_2(\psi)$ in terms of
`approximating' planes. For $\vv a \in \Z^2$ and $ p \in \Z$, let
\begin{equation}
\label{resset4}
 R_{{\vv a},{p}} := \{ {\vv x}  \in \R^2: \
{\vv a}^2\cdot\vv x= p^2  \, \} \ \ .
\end{equation}
It is easily verified, that
$$ \vv x \in S_2(\psi) \hspace{10mm}  {\rm if \ and \ only \ if  }  \hspace{10mm}
\vv x \in \Delta\Big( R_{{\vv a},{ p}} \; ,
\textstyle{\frac{\psi(|\vv a|)}{ {\vv a}.{\vv a}} }  \Big) \cap \I^2
$$ for infinitely many vectors $\vv a\in\Z^2$ and $  p \in
\Z$. The proof of Theorem \ref{t4} follows on establishing the
convergent and divergent parts separately. We make use of the fact
that:
$$
\sum_{h=1}^\infty  \ g\!\left(\dfrac{\psi(h)}{h^2}\right)  \times
\ h^2  \ \ \asymp \ \  \sum_{\vv a\in\Z^2\smallsetminus\{\vv 0\}}
\ g\!\left(\dfrac{\psi(|\vv a|)}{|\vv a|^2}\right)  \times \ |\vv
a| \ .
$$

\vspace{1ex}

\noindent{\em The case of convergence. \ }  The assumption that
$\psi$ is monotonic is irrelevant to this case.  The proof follows
on modifying  the argument of \S\ref{t1conv} in the obvious manner
with $n=2$, $m=1$ and with $\Psi(\vv a)/|\vv a| $ replaced by
$\psi(|\vv a|)/|\vv a|^2$.

\vspace{1ex}

\noindent{\em The case of divergence. \ }
 With reference to our general framework \S\ref{gf}, let $k = 2$ and $m = 1$.
 Hence, $l= 1$. Furthermore, let $ J := \{ (\vv a,  p)
\in   \Z^2\setminus\{\vv 0\}  \times \Z : |\vv a | = |a_1| \} $,
$\alpha := (\vv a,  p) \in J $, $R_\alpha := R_{{\vv a},p}$ where
the latter is given by (\ref{resset4}) and $ \Upsilon_\alpha :=
\psi(|\vv a|)/ {\vv a}.{\vv a}  $. Then,
$$ S_2(\psi)  \  \supset \   \widetilde{S}_2(\psi)
 \ :=  \ \La(\Upsilon)  \, \cap \, \I^2 \ \
. $$ It is easily verified that $ |\widetilde{S}_2(\psi)|_2 = 1 $
whenever $ |S_2(\psi)|_2 = 1 $. Thus,  it suffices to consider the
set $\widetilde{S}_2(\psi)$.  Let $ V : =\{ \vv x = (x_1,x_2) \in
\R^2 \  :\   x_2 = 0  \}  $.  Trivially, conditions (i) and (ii)
of Theorem \ref{t3} are satisfied and with $\Omega = \I^2 $ the
divergence case now follows.

\hfill $\Box$

\subsection{Generalizing Theorem \ref{t3} to fractal subsets $X$  of $\R^k$}

On making use of the general Mass Transference Principle
established in \cite[\S6.1]{mtp} and adapting the Slicing lemma in
an appropriate manner, it is possible to generalize Theorem
\ref{t3} to the following `fractal' setup.  With  $k,l$ and $m$ as
in \S\ref{gf}, let $K$ be a compact subset of $\R^l$. Suppose
there exists   a dimension function $h$ and  constants
$0<c_1<1<c_2<\infty$ and $r_0 > 0$  such that
\begin{equation*}\label{g}
c_1\ h(r)  \ \le \  \cH^h(B(x,r)) \ \le \ c_2\ h(r)  \ ,
\end{equation*}
for any ball $B(x,r)$ with $x\in X$ and $r\le r_0$. In the case $h
: r \to r^{\delta}$ for some $\delta > 0$, the above measure
condition on balls implies that $ \dim K = \delta $ and moreover
that $\cH^{\delta} (K)$ is strictly positive and finite. The
simplest example of a fractal set $K$  satisfying these measure
theoretic  properties is the standard middle third Cantor set  --
simple take $h : r \to r^{\delta}$ with $\delta := \log 2 / \log
3$. More sophisticated  examples include  the attractor $K$
arising from a family of contracting self similarity maps of
$\R^l$ satisfying the open set condition \cite{falc,MAT}. Now let
$$ X \ := \ K \times \R^m \ . $$ Thus, $X$ is a subset of $\R^k$
equipped with the product measure $\mu := \cH^h \times | \ \ |_m
$. Note that if $\dim K = \delta $, then $\dim X =  \delta + m $
and furthermore if $\delta < l $, then $X$ is a set of
$k$--dimensional Lebesgue measure zero. Finally, let $B$ be an
arbitrary ball in $X$ and consider the set $B \cap
\Lambda(\Upsilon) $.  Thus, the points of interest are
 restricted to $X$ since $B$ is by definition a subset of $X$.
 In short, it is possible to
establish an analogue of Theorem \ref{t3} which enables us to
transfer full measure theoretic statements
 with respect to the measure $\mu$ on $X$
 to general Hausdorff measure theoretic
statements for $B \cap \Lambda(\Upsilon)$.
%
The details of this   and its many consequences will be the
subject of a forthcoming article.

\section{Appendix:  Proof of Lemma \ref{slicing1} }

On taking $\phi: \R^k \to V^\perp $ to be the orthogonal projection
map in the following statement, one easily deduces  Lemma
\ref{slicing1}.

\begin{mat} \label{him}
Let $l,k\in\N$ such that $l \leq k $ and $f$ and $g:r\mapsto
r^{-l}f(r)$ be dimension functions. Furthermore, let $A\subset\R^k$
and let $\phi: A \to \R^l $ be a Lipschitz map. Then
$$
\int^{*}_{\R^l} {\cal H}^g(A \cap \phi^{-1}\{y\} ) \; d{\cal L}^ly \
\leq \ \alpha(l) \; 2^l \; {\rm Lip}(\phi)^l  \; {\cal H}^f(A)   \ \
.
$$

\end{mat}

{\em Remark. \ } This is essentially Lemma 7.7 in \cite{MAT}. For
the sake of comparison, the notation adopted above is as far as
possible the same as in \cite{MAT}.  Thus, $\int^{*}$ denotes the
upper integral, ${\cal L}^l$ is the $l$--dimensional Lebesgue
measure on $\R^l$, $\alpha(l) := {\cal L}^l \{x \in \R^l : |x|
\leq 1 \} $ is the volume of the $l$--dimensional unit ball and
${\rm Lip}(\phi)$ is the Lipschitz constant of $\phi$.  To avoid
unnecessary confusion when comparing Lemma \ref{slicing1}* with
Lemma 7.7 in \cite{MAT}, it is worth pointing out that our
statement contains an extra factor of $2^l$ since we have defined
Hausdorff measure in terms of radii of balls rather than
diameters. This extra factor has no effect in deducing Lemma
\ref{slicing1} since all that we require is that the right hand
side of the inequality appearing in  Lemma \ref{slicing1}* is
finite whenever ${\cal H}^f(A) $ is finite.

\vspace*{2ex}

\subsection{Proof of Lemma \ref{slicing1}*}
The statement of Lemma \ref{slicing1}* follows on making the
obvious modifications to the proof of Lemma 7.7 in \cite{MAT}. It
follows from the definition of Hausdorff $f$--measure that for
each $n \in \N$, there exists a cover of  $A$ by  closed balls
$B_{n,1}, B_{n,2}, \ldots $ such that $r(B_{n,i}) \leq 1/n $  and
\begin{equation} \label{a1}
\sum_i f(r(B_{n,i}))     \ \leq \ {\cal H}^{f}_{1/n} (A) \, + \, 1/n
\ \ .
\end{equation}
Let, $F_{n,i} :=  \{ y \in \R^l : B_{n,i} \ \cap  \ \phi^{-1} \{y\}
\neq \emptyset  \}  $. By definition, if $y,z \in F_{n,i}$ then
there exist $u,v \in A \cap B_{n,i} $ such that $\phi(u)=y$ and
$\psi(v)=z$.  It follows that $|y-z| \leq {\rm Lip}(\phi) \, |u-v|$
and so
\begin{equation} \label{a2}
{\cal L}^l(F_{n,i} )  \ \leq \ \alpha(l) \  ({\rm Lip}(\phi) \, 2 \,
r(B_{n,i}) \,  )^l \ .
\end{equation}
For $y \in \R^l$, let
$B_{n,i}^{\phi} (y) $ denote a ball of diameter $d(B_{n,i} \cap
\phi^{-1} \{y\})$ such that $B_{n,i} \cap \phi^{-1} \{y\} \subseteq
B_{n,i}^{\phi} (y) $. On applying Fatou's lemma and using the fact
that   $g $  is non-decreasing, we obtain that
\begin{eqnarray*}
\int^{*}_{\R^l} {\cal H}^g(A \cap \phi^{-1}\{y\} ) \; d{\cal L}^ly
& = &  \int^{*}_{\R^l} \lim_{n \to \infty} {\cal H}^g_{1/n} (A
\cap \phi^{-1}\{y\} ) \; d{\cal L}^ly ~ \hspace*{18ex} ~ \\ & & \\
& \leq & \int_{\R^l} \liminf_{n \to \infty} \sum_i g(\,
r(B_{n,i}^{\phi} (y)) \, ) \; d{\cal L}^ly  \\ & & \\ & \leq &
\liminf_{n \to \infty} \ \sum_i  \ \int_{F_{n,i}}  \!\!\!\! g\Big(
\, \mbox{{\small $\frac{1}{2}$}} \, d(B_{n,i} \cap \phi^{-1}
\{y\}) ) \, \Big) \; d{\cal L}^ly \\ & & \\ & \leq & \liminf_{n
\to \infty} \ \sum_i  \ g(
r(B_{n,i} ) \, ) \  \ {\cal L}^l( F_{n,i} ) \\
& & \\ & \stackrel{(\ref{a2})}{\ \leq \ } & \alpha(l) \  ( 2 \, {\rm
Lip}(\phi) \, )^l \ \liminf_{n \to \infty} \ \sum_i \ f(
r(B_{n,i} ) \, )  \\
& & \\ & \stackrel{(\ref{a1})}{\ \leq \ } & \alpha(l) \  ( 2 \, {\rm
Lip}(\phi) \, )^l \ \liminf_{n \to \infty} \ \Big( {\cal
H}^{f}_{\mbox{{\tiny $1/n$}}}
(A)  +  1/n \Big) \\
& & \\ & \leq & \alpha(l) \  ( 2 \, {\rm Lip}(\phi) \, )^l \ {\cal
H}^{f} (A) \  \ .
\end{eqnarray*}
\hfill $\Box$


\vspace{6mm}

\noindent{\bf Acknowledgments: \, } SV would like to thank the
 Ayesha and Iona for keeping him well  focused on those important  things
 in life -- namely good times, mangoes and  simplicity. Also,
 many thanks to Bridget for sharing nearly half of her years with
 me -- poor thing!

\vspace{5mm}

{ \small

\noindent Victor V. Beresnevich: Department of Mathematics,
University of York,

\vspace{-2mm}

\noindent\phantom{Victor V. Beresnevich: }Heslington, York, YO10
5DD, England.


\noindent\phantom{Victor V. Beresnevich: }e-mail: vb8@york.ac.uk

\vspace{5mm}

\noindent Sanju L. Velani: Department of Mathematics, University
of York,

\vspace{-2mm}

 ~ \hspace{17mm}  Heslington, York, YO10 5DD, England.


 ~ \hspace{17mm} e-mail: slv3@york.ac.uk

 }

\end{document}